\documentclass{jhrs}
\usepackage{amsmath,latexsym,amssymb}

\input xy
\xyoption{all}

\newtheorem{definition}{Definition}[section]
\newtheorem{example}[definition]{Example}
\newtheorem{theorem}[definition]{Theorem}
\newtheorem{proposition}[definition]{Proposition}
\newtheorem{lemma}[definition]{Lemma}
\newtheorem{remark}[definition]{Remark}
\newtheorem{corollary}[definition]{Corollary}
\newtheorem{condition}[definition]{Condition}

\begin{document}
\title{Relative homological categories}
\author{Tamar Janelidze}
\email{tjanelid@maths.uct.ac.za}
\address{Department of Mathematics and Applied Mathematics\\
         University of Cape Town\\
         Rondebosch 7701\\
         Cape Town\\
         South Africa\\}



\begin{abstract} We introduce relative homological and weakly homological
categories $(\mathbf{C},\mathbf{E})$, where ``relative'' refers to a
distinguished class $\mathbf{E}$ of normal epimorphisms in
$\mathbf{C}$. It is a generalization of homological categories, but
also protomodular categories can be regarded as examples. We
indicate that the relative versions of various homological lemmas
can be proved in a relative homological category.
\end{abstract}

\received{June 7, 2006}   
\revised{August 7, 2006}    
\published{Month Day, Year}  
\submitted{Aurelio Carboni}

\startpage{1}

\maketitle

\section{Introduction}
F.\ Borceux and D.\ Bourn \cite{[BB]}, call a category $\mathbf{C}$
\textit{homological} if it is pointed, regular, and protomodular (in
the sense of Bourn \cite{[B]}); in fact they claim that such
categories provide the most convenient setting for non-abelian
versions of various \textit{standard homological lemmas}, such as
snake lemma, 3$\times$3-lemma, etc. Taking this viewpoint, one could
still try, however, to introduce a more general setting involving a
distinguished class $\mathbf{E}$ of regular epimorphisms, where the
homological lemmas are expected to hold only for short exact
sequences $K \rightarrow A \rightarrow B$ with $A \rightarrow B$ in
$\mathbf{E}$. In particular, there is no reason to exclude the
trivial case, where $\mathbf{C}$ is an arbitrary category (say,
pointed and with finite limits and finite colimits) and $\mathbf{E}$
the class of all isomorphisms in $\mathbf{C}$. This idea goes back
to N.\ Yoneda \cite{[Y]}, whose \textit{quasi-abelian} categories
can in fact be defined as pairs $(\mathbf{C},\mathbf{E})$, where
$\mathbf{C}$ is an additive category in which the short exact
sequences $K\rightarrow A\rightarrow B$ with $A\rightarrow B$ in
$\mathbf{E}$ have the same properties as all short exact sequences
in an abelian category.\\ \\
The purpose of this paper is to present a new notion of relative
homological and relative weakly homological categories
$(\mathbf{C},\mathbf{E})$, such that whenever $\mathbf{C}$ is a
pointed category with finite limits and cokernels/coequalizers, we
have:
\begin{itemize}
\item $(\mathbf{C}, \textrm{Isomorphisms in }\mathbf{C})$ always is a
relative homological category;
\item $(\mathbf{C}, \textrm{Split epimorphisms in }\mathbf{C})$ is a
relative weakly homological category if and only if $\mathbf{C}$ is
a protomodular category;
\item $(\mathbf{C}, \textrm{Regular epimorphisms in }\mathbf{C})$ is a
relative homological category if and only if it is a  relative
weakly homological category and if and only if $\mathbf{C}$ is a
homological category;
\item $(\mathbf{C}, \textrm{All morphisms in }\mathbf{C})$ is a relative
homological category if and only if $\mathbf{C}$ is a trivial
category;
\item suitable reformulations of various homological lemmas relative to
$\mathbf{E}$ hold in $\mathbf{C}$.
\end{itemize}

\section{Axioms for relative homological categories}

Throughout the paper we assume that $\mathbf{C}$ is a pointed
category with finite limits and cokernels, and $\mathbf{E}$ is a
class of morphisms in $\mathbf{C}$ containing all isomorphisms.

\begin{definition} \label{def1}
A pair $(\mathbf{C},\mathbf{E})$ is said to be a relative
homological category, if it satisfies the following conditions:
\begin{itemize}
\item[(a)]  The class $\mathbf{E}$ is pullback stable;
\item[(b)] Every morphism in $\mathbf{E}$ is a normal epimorphism;
\item[(c)] $\mathbf{E}$-short-five lemma holds, i.e.\ in every
commutative diagram of the form
$$
\xymatrix{ K \ar@{=}[d] \ar[r]^k & A \ar[r]^f \ar[d]_w & B \ar@{=}[d] \\
           K \ar[r]_{k'} & A' \ar[r]_{f'} & B }
$$
with $k=\mathrm{ker}(f)$, $k'=\mathrm{ker}(f')$, and with $f$ and
$f'$ in $\mathbf{E}$, the morphism $w$ is an isomorphism;
\item[(d)] The class $\mathbf{E}$ is closed under composition;
\item[(e)] If $f \in \mathbf{E}$ and $gf \in \mathbf{E}$, then $g
\in \mathbf{E}$;
\item[(f)]  If a morphism $f$ in $\mathbf{C}$ factors as
$f = em$, in which $e \in \mathbf{E}$ and $m$ is a monomorphism,
then there exists a monomorphism $m'$ in $\mathbf{C}$ and a morphism
$e'$ in $\mathbf{E}$, such that $f = m'e'$;
\item[(g)] If in a commutative diagram
$$
\xymatrix{K \ar[r]^k \ar[d]_u & A \ar[r]^f \ar[d]^w & B \ar@{=}[d] \\
          K' \ar[r]_{k'} & A' \ar[r]_{f'} & B}
$$
$k=\mathrm{ker}(f)$, $k'=\mathrm{ker}(f')$, and morphisms $f$, $f'$,
and $u$ are in $\mathbf{E}$, then $w$ also is in $\mathbf{E}$.
\end{itemize}
We will also say that $(\mathbf{C},\mathbf{E})$ is a relative weakly
homological category whenever it satisfies conditions (a)-(e).
\end{definition}

Note that condition \ref{def1}(c) in fact follows from condition
\ref{def1}(g). Indeed: under the assumptions of \ref{def1}(c), the
morphism $w:A \rightarrow A'$ is in $\mathbf{E}$ since $\mathbf{E}$
contains all isomorphisms and condition \ref{def1}(g) holds. Also,
it is a well known fact that in this situation $\mathrm{ker}(w)=0$.
Since every morphism in $\mathbf{E}$ is a normal epimorphism, we
conclude that $w$ is an isomorphism.

Assuming that condition \ref{def1}(b) holds, we can say that the
conditions/axioms used here are much weaker than those used by G.\
Janelidze, L.\ M\'{a}rki, and W.\ Tholen \cite{[JMT]}. However,
various arguments from \cite{[JMT]}, used there in the proof of the
equivalence of the so-called old and new axioms, can be extended to
our context to obtain various reformulations of our conditions. Some
of them are given in this section.
\begin{condition} \label{con1}
\begin{itemize}
\item[(a)] Every morphism in $\mathbf{E}$ is a regular epimorphism;
\item[(b)] If $f \in \mathbf{E}$ then
$\mathrm{coker}(\mathrm{ker}(f)) \in \mathbf{E}$;
\item[(c)] (``Relative Hofmann's axiom'') If if in a commutative diagram
$$
\xymatrix{ A \ar[r]^f \ar[d]_w & B \ar[d]^{v} \\
           A' \ar[r]_{f'}  & B' }
$$
$f$ and $f'$ are in $\mathbf{E}$, $w$ is a monomorphism, $v$ is
normal monomorphism, and $\mathrm{ker}(f') \leq w$, then $w$ is a
normal monomorphism;
\item[(d)] If in a commutative diagram
$$
\xymatrix@R=10pt{ A \ar[rrd]^{e_1} \ar[dd]_f & \\
           && C \\
           B \ar[rru]_{e_2} }
$$
the morphisms $e_1$ and $e_2$ are in $\mathbf{E}$ and
$\mathrm{Ker}(e_1)=\mathrm{Ker}(e_2)$, then there exists a
factorization \mbox{$f=me$}, in which $m$ is a monomorphism and $e$
is in $\mathbf{E}$.
\end{itemize}
\end{condition}

\begin{theorem} \label{theorem1}
\begin{itemize}
\item[(i)] Condition \ref{def1}(b) implies conditions \ref{con1}(a)
and \ref{con1}(b);
\item[(ii)] Conditions \ref{def1}(a), \ref{def1}(c), \ref{con1}(a)
and \ref{con1}(b) imply condition \ref{def1}(b).
\end{itemize}
\end{theorem}

\begin{proof}
(i) is obvious.

(ii): Let $f:A \rightarrow B$ be a regular epimorphism in
$\mathbf{E}$, and let $k=\mathrm{ker}(f)$ and
\mbox{$q=\mathrm{coker}(k)$}; then condition \ref{con1}(b) implies
that $q$ also is in $\mathbf{E}$. To prove that $f$ is a normal
epimorphism, it is sufficient to show that the canonical morphism
$h: \mathrm{Coker}(k) \rightarrow B$ is an isomorphism. For,
consider the commutative diagram
$$
\xymatrix@C=27pt{ S \ar[d]_{\bar{h}} \ar@<0.5ex>[r]^{s_1}
              \ar@<-0.5ex>[r]_{s_2} & A \ar[r]^>>>>>q \ar[dr]_f &
              \mathrm{Coker}(k) \ar[d]^{h} \\
           R \ar@/_/[ru]^{r_1}
              \ar@/_/@<-1ex>[ru]_{r_2} && B }
$$
in which:
\begin{itemize}
\item[-] $(r_1, r_2)$ is the kernel pair of $f$ and $(s_1,s_2)$ is the
kernel pair of $q$; since the class $\mathbf{E}$ is pullback stable
(condition \ref{def1}(a)), the morphisms $r_1$,$r_2$,$s_1$, and
$s_2$ are in $\mathbf{E}$.
\item[-] $\bar{h}:S \rightarrow R$ is induced by $h$.
\end{itemize}
Since there are canonical isomorphisms
$$ \mathrm{Ker}(s_1) \approx \mathrm{Ker}(q) \approx \mathrm{Ker}(f)
                     \approx \mathrm{Ker}(r_1), $$
we can apply condition \ref{def1}(c) to the diagram
$$
\xymatrix@C=27pt{\mathrm{Ker}(s_1)\ar[d]_{\approx} \ar[r]
                     & S \ar[d]_{\bar{h}} \ar[r]^{s_1} & A \ar@{=}[d]\\
                 \mathrm{Ker}(r_1) \ar[r] & R \ar[r]_{r_1} & A}
$$
This makes $\bar{h}$ an isomorphism; since $f$ and $q$ are regular
epimorphisms, the latter implies that $h$ also is an isomorphism.
\end{proof}

\begin{theorem}\label{theorem2}
\begin{itemize}
\item[(i)] Conditions \ref{def1}(a) and \ref{def1}(c) imply condition
\ref{con1}(c);
\item[(ii)]  Condition \ref{def1}(c) implies
condition \ref{con1}(d);
\item[(iii)] Conditions \ref{def1}(b), \ref{con1}(c),
and \ref{con1}(d) imply condition \ref{def1}(c).
\end{itemize}
\end{theorem}
\begin{proof}
(i): According to the assumptions of \ref{con1}(c), consider the
commutative diagram
$$
\xymatrix{ K \ar@{=}[d] \ar[r]^k & A \ar[r]^f \ar[d]_w & B \ar[d]^v \\
           K \ar[r]_{k'} & A' \ar[r]_{f'} & B' }
$$
in which $f$ and $f'$ are in $\mathbf{E}$, $k'=\mathrm{ker}(f')$,
$k$ is a morphism with $wk=k'$, $w$ is a monomorphism, and $v$ is a
normal monomorphism. It is easy to see that $k$ is in fact the
kernel of $f$, and therefore condition \ref{def1}(c) can be applied
to the diagram
$$
\xymatrix{ K \ar@{=}[d] \ar[r]^k & A \ar[r]^f \ar[d]_{\langle w,f
\rangle} & B
             \ar@{=}[d] \\
           K \ar[r]_>>>>>{h'} & A'{\times}_{B'}B
             \ar[r]_<<<<<{{\pi}_2} & B }
$$
where the projection ${\pi}_2$, being the pullback of $f'$ along
$v$, is in $\mathbf{E}$ by condition \ref{def1}(a). It follows that
$\langle w,f \rangle$ is an isomorphism, and so $w$ is the pullback
of $v$ along $f'$. Since normal monomorphisms are pullback stable,
we conclude that $w$ is a normal monomorphism, as desired.

(ii): Since $\mathbf{E}$ contains all isomorphisms, condition
\ref{con1}(d) follows directly from condition \ref{def1}(c).

(iii): We have to show that if in a commutative diagram
$$
\xymatrix{ K \ar@{=}[d] \ar[r]^k & A \ar[r]^f \ar[d]_w & B \ar@{=}[d] \\
           K \ar[r]_{k'} & A' \ar[r]_{f'} & B }
$$
$f$ and $f'$ are in $\mathbf{E}$, and $k$ and $k'$ are their kernels
respectively, then $w$ is an isomorphism.

It is a well known fact, that in the situation above $w$ has zero
kernel and zero cokernel. By conditions \ref{con1}(d) and
\ref{def1}(b), we have a factorization $w=me$ in which $m$ is a
monomorphism and $e$ is a normal epimorphism in $\mathbf{E}$.
Moreover, since $w$ has zero kernel, $e$ is an isomorphism.
Therefore $w$ is a monomorphism, and applying condition
\ref{con1}(c) to the diagram above, we conclude that $w$ is a normal
monomorphism. Since $w$ has zero cokernel, this implies that $w$ is
an isomorphism, as desired.
\end{proof}
Combining these two theorems, we obtain:
\begin{corollary} \label{cor1}
The following conditions are equivalent:
\begin{itemize}
\item[(i)] A pair $(\mathbf{C},\mathbf{E})$ is a relative weakly homological
category;
\item[(ii)] A pair $(\mathbf{C},\mathbf{E})$ satisfies conditions
\ref{def1}(a), \ref{def1}(c), \ref{def1}(d), \ref{def1}(e),
\ref{con1}(a), and \ref{con1}(b);
\item[(iii)] A pair $(\mathbf{C},\mathbf{E})$ satisfies conditions
\ref{def1}(a), \ref{def1}(b), \ref{def1}(d), \ref{def1}(e),
\ref{con1}(c), and \ref{con1}(d).
\end{itemize}
\end{corollary}
\begin{proof}
The implications (i)$\Rightarrow$(ii), (ii)$\Rightarrow$(i),
(i)$\Rightarrow$(iii), and (iii)$\Rightarrow$(i) follow from
\ref{theorem1}(i), \ref{theorem1}(ii), \ref{theorem2}(i)-(ii), and
\ref{theorem2}(iii) respectively.
\end{proof}

\section{Examples: protomodular and homological categories}

\begin{proposition} \label{prop1} The following conditions are equivalent:
\begin{itemize}
\item[(i)] A pair $(\mathbf{C},\mathbf{E})$ in which
$\mathbf{E}$ is the class of all split epimorphisms in $\mathbf{C}$,
is a relative weakly homological category;
\item[(ii)] $\mathbf{C}$ is a protomodular category in the sense of D.\
Bourn \cite{[B]}.
\end{itemize}
\end{proposition}
\begin{proof}
The implication (i)$\Rightarrow$(ii) follows directly from the
definitions.

(ii)$\Rightarrow$(i): The only condition that requires a
verification here is \ref{def1}(b); however it holds by Proposition
3.1.23 of \cite{[BB]}, which asserts that in a pointed protomodular
category with finite limits, a morphism $f$ is a regular epimorphism
if and only if $f=\mathrm{coker}(\mathrm{ker}(f))$.
\end{proof}

\begin{proposition} \label{prop2} If $\mathbf{C}$ has coequalizers of
kernel pairs and $\mathbf{E}$ is the class of all regular
epimorphisms in $\mathbf{C}$, then the following conditions are
equivalent:
\begin{itemize}
\item[(i)] $(\mathbf{C},\mathbf{E})$ is a relative weakly homological
category;
\item[(ii)] $(\mathbf{C},\mathbf{E})$ is a relative homological category;
\item[(iii)] $\mathbf{C}$ is a homological category in
the sense of F.\ Borceux and D.\ Bourn \cite{[BB]}.
\end{itemize}
\end{proposition}
\begin{proof}
(i)$\Rightarrow$(iii): As follows from (i), the class of all regular
epimorphisms in $\mathbf{C}$ is pullback stable. Therefore, since
$\mathbf{C}$ has kernel pairs and their coequalizers, it admits
(regular epi, mono)-factorization system. Furthermore, using the
same arguments as in \cite{[BJ]}, one can show that in this
situation, protomodularity is equivalent to the $\mathbf{E}$-short
five lemma. It follows that $\mathbf{C}$ is a homological category.

(iii)$\Rightarrow$(ii): Let $\mathbf{C}$ be a homological category
and $\mathbf{E}$ be the class of all regular epimorphisms in
$\mathbf{C}$. Then all properties we need for
$(\mathbf{C},\mathbf{E})$ to be a relative homological category, are
proved in \cite{[BB]}.

Since the implication (ii)$\Rightarrow$(i) is trivial, this
completes the proof.
\end{proof}

\begin{example}
Let $(\mathbf{C},\mathbf{E})$ be a relative weakly homological
category and let $(\mathbf{C'},\mathbf{E'})$ be a pair, in which
$\mathbf{C'}$ is a category with finite limits and $\mathbf{E'}$ is
a class of morphisms in $\mathbf{C'}$ satisfying conditions
\ref{def1}(a), \ref{def1}(d), and \ref{def1}(e). If functor
$F:\mathbf{C} \rightarrow \mathbf{C'}$ preserves finite limits, then
the pair $(\mathbf{C},\mathbf{E} \cap F^{-1}(\mathbf{E'}))$, in
which $F^{-1}(\mathbf{E'})$ is the class of all morphisms $e$ from
$\mathbf{E}$ for which $F(e)$ is in $\mathbf{E'}$, is a relative
weakly homological category. In particular we could take
$\mathbf{C'}$ to be an arbitrary category with finite limits and
$\mathbf{E'}=\mathbf{SplitEpi}$ to be the class of all split
epimorphisms in $\mathbf{C'}$. According to the existing literature
(see e.g.\ \cite{[W]}), an important example is provided by the
forgetful functor $F$ from the homological category $\mathbf{C}$ of
topological groups to the category $\mathbf{C'}$ of topological
spaces; the class $F^{-1}(\mathbf{SplitEpi})$ and the corresponding
concept of exactness play a significant role in the cohomology
theory of topological groups. This also applies to the classical
case of profinite groups, where, however,
$F^{-1}(\mathbf{SplitEpi})$ coincides with the class of all normal
epimorphisms, as shown in Section I.1.2 of \cite{[S]}; another such
result is used in \cite{[C.M]}. The results of \cite{[HW]} also
suggest considering the forgetful functor from the category of
topological groups to the category of groups. On the other hand one
can replace topological groups with more general, so-called
protomodular (=semi-abelian), topological algebras, which form a
homological category due to a result of F.\ Borceux and M.\ M.\
Clementino \cite{[BC]}.
\end{example}
Let us also mention the following ``trivial'' examples:
\begin{example} If $\mathbf{C}$ is an abelian category, and $\mathbf{E}$ is
the proper class of epimorphisms in $\mathbf{C}$ in the sense of
relative homological algebra (see e.g.\ Chapter IX in \cite{[M]})
then $(\mathbf{C},\mathbf{E})$ is a relative weakly homological
category.
\end{example}

\begin{example}
A pair $(\mathbf{C},\mathbf{E})$, in which $\mathbf{E}$ is the class
of all isomorphisms in $\mathbf{C}$, always is a relative
homological category.
\end{example}

\begin{example}
A pair $(\mathbf{C},\mathbf{E})$, in which $\mathbf{E}$ is the class
of all morphisms in $\mathbf{C}$, is a relative homological category
if and only if $\mathbf{C}$ is a trivial category.
\end{example}

\section{Remarks on homological lemmas}
\begin{definition} \label{def2}
Let $(\mathbf{C},\mathbf{E})$ be a relative homological category. A
sequence of morphisms
$$
\xymatrix@!{\ldots \ar[r] & A_{i-1} \ar[r]^{f_{i-1}} & A_{i}
\ar[r]^{f_i} & A_{i+1} \ar[r] & \ldots}
$$
is said to be:
\begin{itemize}
\item[(i)] $\mathbf{E}$-exact at $A_i$, if the
morphism $f_{i-1}$ admits a factorization $f_{i-1}=me$, in which $e
\in \mathbf{E}$ and $m=\mathrm{ker}(f_{i})$;
\item[(ii)] an $\mathbf{E}$-exact sequence, if it is $\mathbf{E}$-exact
at $A_i$ for each $i$ (unless the sequence either begins with $A_i$
or ends with $A_i$).
\end{itemize}
\end{definition}

As easily follows from the definition, the sequence
$$
\xymatrix@!{0 \ar[r] & A \ar[r]^f & B \ar[r]^g & C \ar[r] & 0}
$$
is $\mathbf{E}$-exact, if and only if $f=\mathrm{ker}(g)$ and $g \in
\mathbf{E}$. Having this notion of $\mathbf{E}$-exact sequences, we
may consider the relative cases of various homological lemmas from
\cite{[F.B]} and \cite{[BB]}.\\

\textbf{Relative snake lemma.} Using the same arguments as in the
proof of the Theorem~4.4.2 of \cite{[BB]}, we can prove the
following
\begin{lemma}[Relative snake lemma]\label{lemma1} Let
$(\mathbf{C},\mathbf{E})$ be a relative homological category.
Consider the commutative diagram
$$
\xymatrix@!R@R=3.5pc@!C@C=1pc{& 0 \ar[d] & 0 \ar[d] & 0 \ar[d] \\
          & \mathrm{Ker}(u) \ar[r] \ar[d] &
             \mathrm{Ker}(v) \ar[r] \ar[d] &
             \mathrm{Ker}(w) \ar[d]  \\
           & A \ar[r]^f \ar[d]_u &
             B \ar[r]^g \ar[d]_v & C \ar[r] \ar[d]^w & 0 \\
           0 \ar[r] & {A'} \ar[r]^{f'} \ar[d] &
             {B'} \ar[r]^{g'} \ar[d] & {C'} \ar[d] \\
           & \mathrm{Coker}(u) \ar[r] \ar[d]
           & \mathrm{Coker}(v) \ar[r] \ar[d] & \mathrm{Coker}(w) \ar[d]\\
           & 0 & 0 & 0 }
$$
in which all columns, and the second and the third rows are
$\mathbf{E}$-exact sequences. Suppose $g' \in \mathbf{E}$, or, more
generally, the following conditions hold:
\begin{itemize}
\item[(a)] If ${g_1}' = \mathrm{coker}(f')$ and ${g_2}'$ is the unique morphism with
${g_2}'{g_1}' = g'$, then ${g_1}'$ is in $\mathbf{E}$ and ${g_2}'$
is a monomorphism;
\item[(b)] $(\mathrm{coker}(w))g'_2$ is in $\mathbf{E}$.
\end{itemize}
Then there exists a natural morphism $d:\mathrm{Ker}(w) \rightarrow
\mathrm{Coker}(u)$, such that the sequence
$$
\xymatrix{\mathrm{Ker}(u) \ar[r] & \mathrm{Ker}(v) \ar[r] &
\mathrm{Ker}(w) \ar[r]^d & \mathrm{Coker}(u) \ar[r] &
\mathrm{Coker}(v) \ar[r] & \mathrm{Coker}(w)}
$$
is $\mathbf{E}$-exact.
\end{lemma}

Note that the relative snake lemma can be proved in a relative
weakly homological category $(\mathbf{C},\mathbf{E})$ under some
additional conditions. These additional conditions, however, easily
follow from \ref{def1}(f) and \ref{def1}(g) when
$(\mathbf{C},\mathbf{E})$ is a relative homological category. Being
more precise, to construct the morphism $d$ we only need condition
\ref{lemma1}(a), but to prove the $\mathbf{E}$-exactness of the
sequence above, we need the following:
\begin{itemize}
\item[-] Condition \ref{lemma1}(b);
\item[-] The morphisms $(\mathrm{coker}(v))f'$ and
$(\mathrm{coker}(w))g'_2$ admit a (normal epi in $\mathbf{E}$,
mono)-factorization;
\item[-] If $({B{\times}_{C}(\mathrm{Ker}(w)),{\pi}_1,{\pi}_2})$
is a pullback of $g$ and $\mathrm{ker}(w)$,
\\$(A'{\times}_{(\mathrm{Coker}(u))}(\mathrm{Ker}(w)),{{\pi}_1}',{{\pi}_2}')$
is a pullback of $\mathrm{coker}(u)$ and $d$, and if
$\varphi:{B{\times}_{C}(\mathrm{Ker}(w))} \rightarrow A'$ is induced
by $f'$ and $v{{\pi}_1}:{B{\times}_{C}(\mathrm{Ker}(w))} \rightarrow
B'$, then the morphism ${\langle \varphi,{\pi}_2
\rangle}:B{\times}_{C}(\mathrm{Ker}(w)) \rightarrow
A'{\times}_{(\mathrm{Coker}(u))}(\mathrm{Ker}(w))$ is in
$\mathbf{E}$.\\
\end{itemize}

\textbf{Relative 3$\times$3-lemma.} Using the same arguments as in
the proof of the Theorem~4.5 of \cite{[F.B]}, we can prove the
following
\begin{lemma}[Relative 3$\times$3-lemma]\label{lemma2} Let
$(\mathbf{C},\mathbf{E})$ be a relative homological category. If in
a commutative diagram
$$
\xymatrix{& 0 \ar[d] & 0 \ar[d] & 0 \ar[d] \\
          0 \ar[r] & A \ar[r]^{f} \ar[d]_{u}
          & B \ar[r]^{g} \ar[d]_{v}
          & C \ar[d]^{w} \ar[r] & 0 \\
          0 \ar[r] & A' \ar[r]_{f'} \ar[d]_{u'}
          & B' \ar[r]_{g'} \ar[d]_{v'}
          & C' \ar[d]^{w'} \ar[r] & 0 \\
          0 \ar[r] & A'' \ar[r]_{f''} \ar[d]
          & B'' \ar[r]_{g''} \ar[d]
          & C'' \ar[r] \ar[d] & 0\\
            & 0 & 0 & 0}
$$
the three columns and the middle row are $\mathbf{E}$-exact
sequences, then the first row is an $\mathbf{E}$-exact sequence if
and only if the last row is an $\mathbf{E}$-exact sequence.
\end{lemma}

More generally, the relative version of the 3$\times$3-lemma can be
proved in a relative weakly homological category, if the following
condition holds:
\begin{itemize}
\item[-] The morphism ${\langle v',g' \rangle}:
B' \rightarrow {B''{\times}_{C''}C'}$ (where
$(B''{\times}_{C''}C',{\pi}_1,{\pi}_2)$ is a pullback of $g''$ and
$w'$) is in $\mathbf{E}$.
\end{itemize}
It is easy to show that if $(\mathbf{C},\mathbf{E})$ is a relative
homological category, then this condition follows from Definition
\ref{def1}(g).

\begin{remark}
If $\mathbf{E}$ is the class of all regular epimorphisms in
$\mathbf{C}$, then, by Proposition~\ref{prop2}, $\mathbf{C}$ is a
homological category and all the extra conditions needed for the
proofs of the relative homological lemmas are satisfied. In this
case, the relative snake lemma coincides with its ``absolute
version" proved in \cite{[BB]}, and the relative 3$\times$3-lemma
coincides with its ``absolute version" proved in \cite{[F.B]}.
\end{remark}

\end{document}